\documentclass[10pt]{amsart}
\usepackage{graphicx}%
\usepackage{multirow}%
\usepackage{amsmath,amssymb,amsfonts}%
\usepackage{amsthm}
\usepackage{mathrsfs}%
\usepackage[title]{appendix}%
\usepackage{xcolor}%
\usepackage{textcomp}%
\usepackage{manyfoot}%
\usepackage{booktabs}%
\usepackage{listings}
\usepackage{amssymb}
\usepackage{amsfonts}
\usepackage{amsmath}
\usepackage{epsfig}
\usepackage{indentfirst}
\usepackage[usenames,dvipsnames]{pstricks}
\usepackage{pst-grad}
\usepackage{pst-plot}
\usepackage[margin=1.0in]{geometry}
\usepackage{graphicx}
\usepackage{tikz}
\usepackage{tikz-cd}
\usepackage{float}
\usepackage{tabularx}
\usepackage{mathtools}
\usepackage{bm}
\usepackage{ytableau}

\theoremstyle{thmstyleone}%
\newtheorem{theorem}{Theorem}[section]
\newtheorem{proposition}[theorem]{Proposition}%
\newtheorem{lemma}[theorem]{Lemma}

\theoremstyle{thmstyletwo}%
\newtheorem{example}{Example}%

\theoremstyle{thmstylethree}%
\newtheorem{definition}[theorem]{Definition}%

\DeclareMathOperator {\End}{End}

\DeclareMathOperator{\bc}{\textbf{c}}

\newcommand{\Z}{\mathbb{Z}}
\newcommand{\C}{\mathbb{C}}

\newcommand{\tph}{\tilde{\phi}}
\usetikzlibrary{%
  matrix,%
  calc,%
  arrows%
}

\title{}
\author{}

\begin{document}

\title[Hook fusion procedure for hyper-octahedral groups]{Hook fusion procedure for hyper-octahedral groups}
\author{Dimpi KM}\address{Indian Institute of Science Education and Research Thiruvananthapuram, Email: dtyagi20@iisertvm.ac.in}
\author{Geetha Thangavelu}\address{Indian Institute of Science Education and Research Thiruvananthapuram, Email:  tgeetha@iisertvm.ac.in }


\begin{abstract}
We derive a new expression for the diagonal matrix elements of irreducible representations of the hyper-octahedral group. This expression is obtained using Grime's hook fusion procedure for symmetric groups, which minimizes the number of auxiliary parameters required in the fusion process.
\end{abstract}


\maketitle
\section{Introduction}\label{sec1}

The study of the Yang-Baxter equation led to the discovery of quantum groups, which have numerous applications in mathematics and mathematical physics. Given a fundamental solution to the Yang-Baxter equation, the fusion procedure enables the construction of new solutions, known as \textit{fused solutions}. Starting from an initial solution on the space $V$, the fusion procedure generates fused solutions on tensor products of $V$ and its subspaces.

Let $V$ be a finite-dimensional complex vector space. Let $R(u)$ be a function of a single variable $u\in \C$ with values in $\End(V\otimes V)$, where $u$ is the spectral parameter. The function $R$ is called a solution of the Yang-Baxter equation (YBE) if it satisfies the following relation:
\begin{equation}
R_{12}(u)R_{13}(u+v)R_{23}(v) = R_{23}(v)R_{13}(u+v)R_{12}(u)\label{YBE},
\end{equation}
where both sides of the equation belong to $\End(V^{\otimes 3})$, and the subscripts of $R(u)$ indicate the copies of $V$ on which the operator $R(u)$ acts. For example, $R_{12}(u):=R(u)\otimes \operatorname{Id}_V$ means that $R(u)$ acts on the first two factors and $\operatorname{Id}_V$ acts on the third. In this context, by a slight abuse of notation we may refer to the operator $R(u)$ itself as a solution of the Yang-Baxter equation. The simplest solution of (\ref{YBE}), known as the Yang $R$-matrix, is given by
\begin{equation}
R(u)=\operatorname{Id}_{V^{\otimes 2}}-\frac{P}{u}.\label{RM},
\end{equation}
where the operator $P$ is the permutation operator $P(x\otimes y) = y\otimes x$ for $x, y\in V$.

Let $n\geq 2$ be an integer, and let $\bc := (c_1,\dots,c_n)$ be an $n$-tuple of complex parameters. The {\it fusion function}, denoted by $F(\bc)$, is an endomorphism of $V^{\otimes n}$ defined by
$$
F(\bc):=\prod_{1\leq i<j\leq n}^{\rightarrow}R_{ij}(c_i-c_j),
$$
where the factors are taken in lexicographical order. The function $R(u)$, and hence the fusion function $F(\bc)$, may have singularities. Choose $\bc$ such that $F(\bc)$ is well-defined. The image $W_{\bc}$ of $F(\bc)$ is a subspace of $V^{\otimes n}$.  From an arbitrary solution $R(u)$ of the Yang-Baxter equation on $V\otimes V$, a \textit{fused solution} $R_{\bc,\underline{\bc}}(u)$ is constructed using two sequences of complex numbers $\bc := (c_1,\dots, c_n)$ and ${\underline{\bc}} := (c_{1'},\dots, c_{n'})$, for more details see \cite{PL17}. The resulting solution $R_{\bc,\underline{\bc}}(u)$ is an endomorphism of the vector space $V^{\otimes n}\otimes V^{\otimes n'}$. It has been shown that the fused operator $R_{\bc,\underline{\bc}}(u)$ preserves the subspace $W_{\bc}\otimes W_{\underline{\bc}}$ of $V^{\otimes n}\otimes V^{\otimes n'}$, where 
\begin{equation}
W_{\textbf{c}} := (\prod_{1\leq i< j\leq n}^{\rightarrow} R_{ij}(c_i -c_j))V^{\otimes n}\label{W}.
\end{equation}
For the solution of the form (\ref{RM}), the expression for the fusion function becomes $$F(\bc):= \prod_{1\leq i< j\leq n}^{\rightarrow}(\operatorname{Id}_{V^{\otimes n}} - \frac{P_{ij}}{c_i - c_j}),$$ where $P_{ij}$ denotes the permutation operator acting on the $i$-th and $j$-th factors.   

On the other hand, the symmetric group $S_n$ acts on $V^{\otimes n}$ by permuting the tensor factors. The fusion function $F(\bc)$ coincides with the image in $\End(V^{\otimes n})$ of the following rational function with values in $\C[S_n]$:
$$
\prod_{1\leq i<j\leq n}^{\rightarrow}\phi_{ij}(c_i,c_j),
$$
where $\phi_{ij}(u,v)=1-\frac{(ij)}{u-v}$. This rational function is referred to as the ``fusion function'' for the symmetric groups.  A complete set of pairwise orthogonal primitive idempotents can be constructed by taking the limiting values of this rational function. Specifically, the values $c_i$ are chosen to be the respective \textit{contents} of the standard tableau $\bm{T}$ associated to a partition of $n$. 

 This approach provides an expression for the diagonal matrix elements of the irreducible representations of $S_n$. It originates from the work of Jucys \cite{Ju66}, was further developed by Cherednik \cite{Che86}, and is known as the \textit{fusion procedure} for symmetric groups. The complete proof was later given by Nazarov \cite{Naz97}. Further insights into this procedure can be found in various works, including the work of Molev \cite{Mo08}, who introduced an alternative approach using Jucys-Murphy elements, offering a relatively simple proof of the fusion procedure for the symmetric groups. The fusion procedure has been extended to various other structures, such as the Hecke algebra of type $A$, the Brauer algebra, and the BMW algebra \cite{Gri07,IM08,IM10,IM14}.

The fusion procedure requires $n$ auxiliary parameters. Nazarov introduced a distinct limiting procedure that reduces the number of auxiliary parameters. Later, Grime \cite{Gri05} developed a version of the fusion procedure called the \textit{hook fusion procedure}, which leads to a further reduction in the number of auxiliary parameters.

In this article, we adopt Grime's approach to derive an analogue of the hook fusion procedure for the hyper-octahedral groups. We derive the hook fusion procedure for the Coxeter groups of type $B$, together with a description of a complete set of pairwise orthogonal primitive idempotents for $B_n$. While the standard fusion procedure for the hyper-octahedral groups and complex reflection groups of type $G(m,1,n)$ was developed in \cite{OP14}, our approach minimizes the number of auxiliary parameters and simplifies the construction of idempotents.

\section{Preliminaries}
Let us fix some notations and recall well-known facts about the representations of the hyperoctahedral groups; see \cite{Pu99} for more details.
\begin{definition}
The hyperoctahedral group $\Z_{2}\wr S_n$ is a Coxeter group of type B, generated by the elements $s_{1}, s_{2},\dots, s_{n-1}$ and $t$ with relations:  

\begin{equation}
\begin{split}
\begin{aligned}
&s_i^{2} = 1, \quad t^{2} = 1 && \text{for $i = 1,\dots, n-1$}, \\
&s_{i}s_{i+1}s_{i} = s_{i+1}s_{i}s_{i+1} &&\text{for $i = 1,\dots, n-2$}, \\
&s_{i}s_{j} = s_{j}s_{i} &&\text{for $i, j = 1,\dots, n-1$ such that $|i-j|>1$}, \\ 
&s_{i}t = ts_{i} && \text{for $i = 2,\dots, n-1$},\\
&ts_{1}ts_{1} = s_{1}ts_{1}t.   
\end{aligned} \label{re}
\end{split} 
\end{equation}

\end{definition}

We set:
\begin{align}
t(p):= &\frac{1}{2}(1+pt), \label{tre}
\end{align}
and 
\begin{align}
S_{i}(a,a',p,p') := s_{i} + \frac{\delta_{p,p'}}{a-a'}.  \label{S}
\end{align}
Here, $\delta_{p,p'}$ denotes the Kronecker delta, and for $p=p'$, the elements in equation (\ref{S}) correspond to the analogs of \textit{Baxterized elements} for the hyperoctahedral group.

A \textit{partition} of an integer $n$, denoted by $\lambda$, is a sequence of positive integers $\lambda:=(\lambda_1,\lambda_2,\dots,\lambda_k)$ such that $\lambda_1\geq \lambda_2\geq \dots \geq \lambda_k$ and $\lambda_1+\lambda_2+\dots+\lambda_k=n$, which is represented as $\lambda\vdash n.$ This definition is extended to a \textit{bi-partition} of $n$, which consists of a pair of partitions $\bm{\lambda}=(\lambda^1,\lambda^2)$, where $\lambda^1=(\lambda^1_1,\lambda^1_2,\dots,\lambda^1_k)\vdash m$ and $\lambda^2=(\lambda^2_1,\lambda^2_2,\dots,\lambda^2_t)\vdash l$, with $m+l=n$.

Similarly, for a bi-partition $\bm{\lambda}=(\lambda^1,\lambda^2)\vdash n$, we define two Young diagrams, one for each of the partitions $\lambda^1$ and $\lambda^2$, respectively. This pair of diagrams we call the \textit{Young diagram} associated with the bi-partition $\bm{\lambda}$. Hereafter, the term ``Young diagram" will refer to the associated Young diagram of the given partition or bi-partition, depending on the context.

Given a partition $\lambda:=(\lambda_1,\lambda_2,\dots,\lambda_k)\vdash n$, the associated \textit{Young diagram} consists of $k$ rows of boxes, where $j$-th row contains $\lambda_j$ left-justified nodes, where a \textit{node} $\alpha$ in a Young diagram is a specific box in the array and can be described by the co-ordinates $(i,j)$. The content of a node $\alpha$ in a Young diagram is given by the formula:
$$c(\alpha):=j-i.$$

In the context of Young diagrams of a bi-partition of shape $\bm{\lambda}$, a node of a bi-partition denoted by $\bm{\alpha}$, is a pair $(\alpha , k)$, where $\alpha$ is a node of the Young diagram of $\lambda^k$, and $k\in\{1,2\}$ refers to the position of the node in the Young diagram of a bi-partition.  Again define $c(\bm\alpha):=c(\alpha)=j-i$ to be the \textit{content} of the node $\alpha$ in a Young diagram of a bi-partition. Let $\{\xi_1,\xi_2\}$ be the set of distinct roots of unity, ordered arbitrarily. We define $p(\bm{\alpha}):=\xi_k.$ The notion of removing or adding a node in a Young diagram extends naturally to the Young diagrams of a bi-partition.

Given a partition $\lambda\vdash n$, a \textit{standard tableau} $T$ is a filling of the nodes of its Young diagram with the numbers $1,2,\dots,n$, each occurring exactly once, such that the entries in each row and column are strictly increasing. Similarly, a \textit{standard bi-tableau} $\bm{T}=(T_1,T_2)$ of shape $\bm{\lambda}=(\lambda^1,\lambda^2)$ is obtained by filling the nodes of the two Young diagrams collectively with the numbers $1,2,\dots,n$, each occurring exactly once, such that the entries increase along each row and column in each component. Let $c_i$ denote the content of the node containing $i$ in $\bm{T}$, and let $p_i=\xi_k$ if $i$ appears in the $k$-th component.

In the Young diagram of shape $\lambda$, the $(i,j)$-\textit{hook} of a node $(i,j)$ is the set of boxes in $\lambda$
$$\{(i,j'):j'\geq j\}\cup\{(i',j):i'\geq i\}.$$ That is, the set of nodes that are in the same row as that of $(i,j)$ and to the right of it, along with the nodes that are in the same column as $(i,j)$ and below it. We call the $(i,i)$-hook the $i$-th \textit{principal hook}. The \textit{hook length} of a node $\alpha$, denoted by $h(\alpha)$, is the number of nodes in the hook of $\alpha$. The above definitions extend componentwise to the Young diagrams of bipartitions. For a bipartition $\bm{\lambda}$, we set
\begin{align}
f_{\bm{\lambda}} := \left(
\prod_{\bm{\alpha}\in\bm{\lambda}}
h_{\bm{\lambda}}(\bm{\alpha})
\right)^{-1},
\label{ho}
\end{align}
where $h_{\bm{\lambda}}(\bm{\alpha})$ denotes the hook length of the node $\bm{\alpha}$ in its corresponding component of $\bm{\lambda}$.

The irreducible representations of $\mathbb{Z}_{2}\wr S_{n}$ are indexed by bi-partitions of $n$. For further details, we refer to \cite{Ar97,Ca98,Pu99} and \cite{Ho74}.  The seminormal basis elements of an irreducible representation, labeled by a bi-partition $\bm{\lambda}$, are indexed by standard bi-tableaux of shape $\bm{\lambda}$. Given a standard bi-tableau $\bm{T}$, let $v_{\bm{T}}$ denote the corresponding basis vector in the irreducible representation $V_{\bm{\lambda}}$. The \textit{diagonal matrix element} $F_{\bm{T}}$ associated with $v_{\bm{T}}$ is defined as
\begin{align}
F_{\bm{T}} = \sum_{g\in \mathbb{Z}_{2}\wr S_{n}} \langle v_{\bm T},gv_{\bm T}\rangle g \in \mathbb{C}[\mathbb{Z}_{2}\wr S_{n}].\label{d}
\end{align}
where $< v_{\bm{T}}, gv_{\bm{T}}>$ is the matrix coefficient.

The primitive idempotent corresponding to $v_{\bm{T}}$ is the following element in the group algebra $ \C[\Z_{2}\wr S_{n}]$:
\begin{align}
E_{\bm{T}} = \frac{d_{\bm{\lambda}}}{|\mathbb{Z}_{2}\wr S_{n}|}\sum_{g\in \mathbb{Z}_{2}\wr S_{n}} \langle v_{\bm T},gv_{\bm T}\rangle g \in \C[\Z_{2}\wr S_{n}]\label{mu}
\end{align}
where $d_{\bm{\lambda}}=\dim(V_{\bm{\lambda}})$ is the dimension of the irreducible representation $V_{\bm{\lambda}}$.

Consider a standard bi-tableau $\bm{T} = (T_{1},T_{2})$, and let $\bm{T}'=s_{i}\bm{T}$, where $\bm{T}'$ is obtained by exchanging the entries $i$ and $i+1$ in $\bm{T}$. If $c_{i}(\bm{T})\neq c_{i+1}(\bm{T})$, we define 
\begin{align}
d_{i}(\bm{T})=(c_{i+1}(\bm{T})-c_{i}(\bm{T}))^{-1}. 
\end{align}
When $\bm{T}'$ is non-standard, it means that $i$ and $i+1$ are adjacent either in the same row or column of $T_1$ or $T_2$, and in this case $d_{i}(\bm{T})=1$ or $d_{i}(\bm{T})=-1$. Following Young's construction, the vectors of the Young basis can be normalised so that the action of the generators of $\C[\Z_{2}\wr S_{n}]$ on the semi-normal basis is as follows:\\
(i) If $i$ and $i+1$ lie in the same tableau, then
\begin{equation}
 s_i.v_{\bm{T}} =
  \begin{cases}
    d_{i}(\bm{T})v_{\bm{T}} + \sqrt{1-d_{i}(\bm{T})^2}v_{\bm{T}'} & \text{if $\bm{T}'$ is standard }, \\
    d_{i}(\bm{T})v_{\bm{T}} &\text{if $\bm{T}'$ is non-standard}.\label{A1}
  \end{cases}
\end{equation}
(ii) If $i$ and $i+1$ lie in different tableaux, then
\begin{equation}
s_{i}.v_{\bm{T}} = v_{\bm{T}'}.\label{A2}
\end{equation}
(iii) The action of $t$ is given by 
\begin{equation}
 t.v_{\bm{T}} =
  \begin{cases}
    v_{\bm{T}} & \text{if 1 lies in $T_1$}, \\
    -v_{\bm{T}} & \text{if 1 lies in $T_2$}.\label{A3}
  \end{cases}
\end{equation}\\
Moreover, if $\bm{T}'$ is standard and $i$ and $i+1$ are in the same tableau, then from (\ref{A1}) we also have the relation 
\begin{align}
\sqrt{1-d_{i}(\bm{T})^2}v_{\bm{T}'} = (s_{i}-d_{i}(\bm{T}))v_{\bm{T}}.\label{Re1}
\end{align}
Using the definition of the diagonal matrix element $F_{\bm{T}}$ from (\ref{d}), this leads to the important relation
\begin{align}
(1-d_{i}(\bm{T})^2)F_{\bm{T}'} = (s_{i}-d_{i}(\bm{T}))F_{\bm{T}}(s_{i}-d_{i}(\bm{T})).\label{Re2}
\end{align}
Moreover, if $\bm{T}'$ is standard and $i$ and $i+1$ lie in different tableaux, then, by (\ref{A2}), we obtain the following identity: 
\begin{align}
F_{\bm{T}'} = s_{i}F_{\bm{T}}s_{i}.\label{Re3}
\end{align}
In this article, our aim is to develop a hook fusion procedure for $\mathbb{Z}_{2}\wr S_{n}$, focusing on minimizing the number of auxiliary parameters typically used in the fusion procedure for $\mathbb{Z}_{2}\wr S_{n}$.

To set this up, we define $\phi_{1}(u,v) := t(v)$ and, for $k=2,\dots,n$, we define
\begin{equation}
\begin{split}
\begin{aligned}
\phi_{k}(u_{1},\dots, u_{k-1}, u_{k}, v_{1}, \dots, v_{k-1}, v_{k}):= & S_{k-1}(u_{k}, u_{k-1}, v_{k}, v_{k-1})\phi_{k-1}(u_{1},\dots, u_{k-2}, u_{k}, v_{1}, \dots, v_{k-2}, v_{k})s_{k-1}\\
= & S_{k-1}(u_{k}, u_{k-1}, v_{k}, v_{k-1})\dots S_{1}(u_{k}, u_{1}, v_{k}, v_{1})t(v_{k})s_{1}\dots s_{k-1}. \label{Rf}  
\end{aligned}
\end{split}
\end{equation} 

Using these, we define the rational function with values in $\mathbb{C}[\mathbb{Z}_{2}\wr S_{n}]$: 
\begin{align}
\Phi(u_{1},\dots, u_{n},v_{1},\dots, v_{n}):= \prod_{k=1,\dots,n}^{\leftarrow}\phi_{k}(u_{1},\dots, u_{k}, v_{1}, \dots, v_{k}).\label{Rf1}
\end{align}
where the product is ordered from right to left, as indicated by the arrow.

It was shown in \cite{OP14} that if we consecutively evaluate this rational function (\ref{Rf1}) as
\begin{align}
\Phi(u_{1},\dots, u_{n},v_{1}, \dots, v_{n})|_{v_{i}=p_{i}, i = 1,\dots, n}|_{u_{1}=c_{1}}\dots |_{u_{n}=c_{n}} \label{Rf'},
\end{align} 
then this coincides with $\frac{1}{f_{\bm{\lambda}}}E_{\bm{T}}$, where $\bm{T}$ is a standard bi-tableau of shape $\bm{\lambda}$, and $p_i$ and $c_i$ are defined earlier.

To compute $E_{\bm T}$ in an alternative way, we first rewrite the rational function (\ref{Rf1}) as:
\begin{align}
\Phi_{\bm{T}}(z_{1},\dots, z_{n}, w_{1},\dots, w_{n}):= \prod_{k=1,\dots,n}^{\leftarrow}\phi_{k}(c_{1}+z_{1},\dots, c_{k}+z_{k},p_{1}+w_{1}\dots, p_{k}+w_{k}). \label{RF1}
\end{align}
Note that if the entries $k$ and $l$ $(k<l)$ lie on the same diagonal in the same component of $\bm{T}$, then $\Phi_{\bm{T}}(z_{1},\dots, z_{n}, w_{1}\dots, w_{n})$ has a singularity at $z_{k} = z_{l}$.

To address this, we define the set
$$
H_{\bm{T}} :=\left\{\begin{aligned}
(z_{1},\dots, z_{n}) \in \C^{n} \;\big|\;& z_{i} = z_{j} \quad \text{whenever $i$ and $j$ lie in the same principal hook}\\ &\text{of a component of $\bm{T}$}
\end{aligned}\right\}.
$$
where the \textit{principal hook} refers to the hook associated with a diagonal node.

We will show that when $\Phi_{\bm{T}}(z_{1},\dots, z_{n}, w_{1}\dots, w_{n})$ is restricted to $H_{\bm{\lambda}}$, it becomes regular at $z_{1} = z_{2} = \dots = z_{n}$, and furthermore, its value at $z_{1} = z_{2} = \dots = z_{n}$ and $w_{1} = w_{2} = \dots = w_{n} = 0$ is a scalar multiple of a primitive idempotent of $\mathbb{C}[\mathbb{Z}_{2}\wr S_{n}]$. The following is the main result of this article, which we will prove in the next section.

\begin{theorem}
Let $\bm{\lambda}$ be a bi-partition of $n$, and let $\bm{T}$ be a standard bi-tableau of shape $\bm{\lambda}$. Then the restriction of the rational function $\Phi_{\bm{T}}(z_1,\dots,z_n,w_1,\dots,w_n)$ to $H_{\bm{T}}$ is regular at $z_1=z_2=\dots=z_n$. Moreover, the value of this restriction at $z_1=z_2=\dots=z_n$ and $w_1=w_2=\dots=w_n=0$ coincides with $\frac{1}{2^n}F_{\bm{T}}$.
\end{theorem} 

\section{Hook Fusion Procedure For $\mathbb{C}[\mathbb{Z}_{2}\wr S_{n}]$}

\begin{definition}\label{def:hook-bitableau}
Let $\bm{\lambda}=(\lambda^1,\lambda^2)$ be a bi-partition of $n$. We construct a bi-tableau as follows. In the first diagram $\lambda^1$, fill the column of the first principal hook with the first consecutive integers and then fill the remaining boxes in its row with the next consecutive integers. Repeat this procedure for the second principal hook, and continue through all the principal hooks of $\lambda^1$. Then apply the same procedure to the diagram $\lambda^2$, using the remaining integers. The resulting standard bi-tableau, denoted by $\bm{T}^{0}$, is called the \textit{hook bi-tableau} of shape $\bm{\lambda}$.
\end{definition}

\begin{example}
The hook bi-tableau of shape $\bm{\lambda}=((3,2,1),(3,3,2))$ is
$$
\left(\begin{ytableau}
       1&4&5\\
       2&6\\
       3
\end{ytableau},\quad
\begin{ytableau}
       7&10&11\\
       8&12&14\\
       9&13
\end{ytableau}\right).
$$
\end{example}

In this section, our aim is to prove Theorem 2.2. To do this, rather than working with an arbitrary standard bi-tableau $\bm{T}$ of shape $\bm{\lambda}$, we will focus on the hook bi-tableau $\bm{T}^{0}$ of shape $\bm{\lambda}$. We will show that the regularity of the rational function associated with a standard bi-tableau $\bm{T}$ of shape $\bm{\lambda}$ is equivalent to the regularity of the same function for the hook bi-tableau $\bm{T}^{0}$ of shape $\bm{\lambda}$.

We shall use the following relations for the Baxterized elements:
\begin{align}
S_{i}(a,a',p,p')S_{i+1}(a,a'',p,p'')S_{i}(a',a'',p',p'') = S_{i+1}(a',a'',p',p'')S_{i}(a,a'',p,p'')S_{i+1}(a,a',p,p') \label{Bx1}
\end{align}
and
\begin{align}
S_i(a,a',p,p')S_j(b,b',q,q') =S_j(b,b',q,q')S_i(a,a',p,p') \label{Bx2}
\end{align}
for $|i-j|>1$, and
\begin{align}
S_i(a,a',p,p')S_i(a',a,p',p)
=1-\frac{\delta_{p,p'}}{(a-a')^2}.\label{Bx3}
\end{align}

For convenience, we introduce the notation:$$S_{k}(z_{j}+c_{j}, z_{i}+c_{i}) := S_{k}(z_{j}+c_{j}, z_{i}+c_{i}, w_{j}+p_{j}, w_{i}+p_{i}).$$ Note that this factor $S_{k}(z_{j}+c_{j}, z_{i}+c_{i})$ is not defined when $z_j = z_i $ and this occurs if and only if the entries $j$ and $i$ lie on the same diagonal within the same tableau of $\bm{T}$. Thus, we refer to such a pair $(j,i)$ as a \textit{singularity} and the term $S_{k}(z_{j}+c_{j}, z_{i}+c_{i})$ as a \textit{singularity term} or simply \textit{singularity}.

Let $i$ and $j$ $(i<j)$ be in the same hook of a bi-tableau $\bm{T}$. If $i$ and $j$ are consecutive in the column of the hook, then on $H_{\bm{T}}$, $S_{i}(z_{j}+c_{j}, z_{i}+c_{i})$ takes the value $(s_i - 1)$. Similarly, when $i$ and $j$ are consecutive in the same row of the hook, $S_{i}(z_{j}+c_{j}, z_{i}+c_{i})$ equals $(s_i + 1)$ on $H_{\bm{T}}$. Consequently, $-\frac{1}{2}S_{i}(z_{j}+c_{j}, z_{i}+c_{i})$ acts as an idempotent when $i$ and $j$ are in the same column next to each other, and $\frac{1}{2}S_{i}(z_{j}+c_{j}, z_{i}+c_{i})$ acts as an idempotent when $i$ and $j$ are in the same row next to each other.

The inclusion of singularity terms in the product $\Phi_{\bm{T}}(z_1,\dots, z_n, w_1,\dots, w_n)$ implies that this product may or may not be well defined on the vector subspace $H_{\bm{T}}$ in the limit where $z_1 = \dots = z_n$. The following lemma will be useful to show that $\Phi_{\bm{T}}(z_1,\dots, z_n, w_1,\dots, w_n)$ is indeed well-defined on the subspace, with all singularities being removable. 

\begin{lemma}\label{Reg 1}
The product
$$
\left(s_{i+1}+\frac{1}{u-v}\right)
\left(s_i+\frac{1}{u-w}\right)
\left(s_{i+1}+\frac{1}{v-w}\right)
$$
is regular at $u=w$ when restricted to $v=w\pm1$.
\end{lemma}

\begin{proof}
Suppose that $v=w\pm1$. Using $s_{i+1}s_is_{i+1}=s_is_{i+1}s_i$, the product can be written as
$$
\left(s_{i+1}s_i+\frac{s_i\pm1}{u-v}\right)(s_{i+1}\pm1).
$$
At $u=w$ we have $u-v=\mp1$, so the last expression is regular.
\end{proof}

Next, we define some special elements of $\Z_{2}\wr S_{n}$ inductively as:
\begin{align}
j_{1} = t,\quad j_{i+1} = s_{i}j_{i}s_{i}  \qquad  \text{ for }   i = 1,\dots, n-1.\label{JM}
\end{align}

These elements $j_{i}$ form a set analogous to the Jucys-Murphy elements for $\Z_{2}\wr S_{n}$. They commute with all the Coxeter generators $s_{k},$ except for $s_{i}$ and $s_{i-1},$ specifically: 
\begin{align}
j_{i}s_{k} = s_{k}j_{i} \qquad if \quad k\neq i, i-1.\label{Jm1}
\end{align}

Furthermore, the following relation holds:
\begin{align}
j_{l}s_{l}s_{l+1}\dots s_{n-1} = s_{l}s_{l+1}\dots s_{n-1}j_{n}.\label{Jm2}
\end{align}

To ensure clarity in later computations, we introduce a generalized version of the expression given in (\ref{tre}),
\begin{align}
J_{i}(p) = \frac{1}{2}(1 + pj_{i}).\label{Jm3}
\end{align}

Observe, 
\begin{align}
J_{i}(p)J_{k}(q) = J_{k}(q)J_{i}(p) \quad for \quad i\neq k.\label{Jm4}
\end{align}

\begin{proposition}\label{prop:factorization}
The rational function $\Phi_{\bm{T}}(z_{1},\dots, z_{n}, w_{1}\dots, w_{n})$ can be expressed as  
\begin{align}\label{RF3}
\Phi_{\bm{T}}(z_{1},\dots, z_{n}, w_{1}\dots, w_{n}) = (\prod_{k=1,\dots,n}^{\leftarrow}&\tph_{k}(z_{1}+c_{1},\dots, z_{k}+c_{k},w_{1}+p_{1}\dots, w_{k}+p_{k}))  \\\nonumber\noindent & J_{1}(w_1+p_1)\dots J_{n}(w_n+p_n) 
\end{align} 
where 
\begin{align}\label{Rf4}
\tph_k(z_{1}+c_{1},\dots, z_{k}+c_{k},w_{1}+p_{1}\dots, w_{k}+p_{k}) = &S_{k-1}(z_{k}+c_k, z_{k-1}+c_{k-1})\dots \\\nonumber\noindent & S_{1}(z_{k}+c_k, z_{1}+c_1)s_{1}\dots s_{k-1}. 
\end{align}
\end{proposition}

\begin{proof}
By definition, $$\Phi_{\bm{T}}(z_{1},\dots, z_{n}, w_{1}\dots, w_{n}) = \prod_{k=1,\dots,n}^{\leftarrow}\phi_{k}(z_{1}+c_{1},\dots, z_{k}+c_{k},w_{1}+p_{1}\dots, w_{k}+p_{k})$$
where,
\begin{align*}
\phi_{k}(z_{1}+c_{1},\dots, z_{k}+c_{k},w_{1}+p_{1},\dots, w_{k}+p_{k}) = & S_{k-1}(z_{k}+c_k, z_{k-1}+c_{k-1})\dots S_{1}(z_{k}+c_k, z_{1}+c_1)\\\nonumber\noindent & t(w_k+p_k)s_{1}\dots s_{k-1}.
\end{align*}
Since $t(w_k+p_k) = J_1(w_k+p_k)$, and using (\ref{Jm2}) we can write
$$\phi_{k}(z_{1}+c_{1},\dots, z_{k}+c_{k},w_{1}+p_{1}\dots, w_{k}+p_{k}) = \tph_{k}(z_{1}+c_{1},\dots, z_{k}+c_{k},w_{1}+p_{1}\dots, w_{k}+p_{k})J_k(w_k + p_k).$$ 
From (\ref{Jm1}), we have $J_k(w_k+p_k)S_i(z_j+c_j, z_i+c_i) = S_i(z_j+c_j, z_i+c_i)J_k(w_k+p_k)$ when $i\neq k, k-1$. Hence, $J_k(w_k+p_k)\tph_{k-1}(z_{1}+c_{1},\dots, z_{k-1}+c_{k-1},w_{1}+p_{1}\dots, w_{k-1}+p_{k-1}) = \tph_{k-1}(z_{1}+c_{1},\dots, z_{k-1}+c_{k-1},w_{1}+p_{1}\dots, w_{k-1}+p_{k-1})J_k(w_k+p_k)$. More generally, if $l<k$, then $\tph_l$ involves only the generators $s_1,\dots,s_{l-1}$, and hence $J_k(w_k+p_k)$ commutes with $\tph_l$. Therefore, each $J_k(w_k+p_k)$ can be moved through all the factors $\tph_{k-1},\dots,\tph_1$. Applying (\ref{Jm4}) we get the relation (\ref{RF3}) and the proof is completed.     
\end{proof}

We now give the proof of the first statement of Theorem 2.2.

\begin{proposition}\label{prop:regularity}
The restriction of the rational function $\Phi_{\bm{T}}(z_{1},\dots, z_{n}, w_{1}\dots, w_{n})$ to the subspace $H_{\bm{T}}$ is regular at $z_{1} = z_{2} = \dots = z_{n}.$
\end{proposition}

\begin{proof}
The proof consists of two steps. We first transfer the regularity question from an arbitrary standard bi-tableau to the hook bi-tableau of the same shape. We then remove the possible singularities in the function associated with the hook bi-tableau.

\medskip
\noindent
\textbf{Step 1: Reduction to the hook bi-tableau.}

Suppose that $\bm{T}'=s_i\bm{T}$ is standard. Thus, $\bm{T}'$ is obtained from $\bm{T}$ by interchanging the entries $i$ and $i+1$. Let
$$
(z'_1,\ldots,z'_n)
\quad\text{and}\quad
(w'_1,\ldots,w'_n)
$$
be obtained from $(z_1,\ldots,z_n)$ and $(w_1,\ldots,w_n)$ by interchanging the variables with indices $i$ and $i+1$.

For brevity, let $\tph_k^{\bm{T}}$ denote the factor $\tph_k$ whose arguments are determined by the contents and position parameters of $\bm{T}$. The braid and commuting relations for the Baxterized elements imply that, for $k>i+1$,
\begin{align}
&S_i(z_i+c_i(\bm{T}),z_{i+1}+c_{i+1}(\bm{T}))
\tph_k^{\bm{T}}(z_1,\ldots,z_k,w_1,\ldots,w_k)
\nonumber\\
&\qquad=
\tph_k^{\bm{T}'}(z'_1,\ldots,z'_k,w'_1,\ldots,w'_k)
S_i(z_i+c_i(\bm{T}),z_{i+1}+c_{i+1}(\bm{T})).
\label{E1}
\end{align}
For the two factors $\tph_{i+1}^{\bm{T}}$ and $\tph_i^{\bm{T}}$, the corresponding relation is
\begin{align}
&S_i(z_i+c_i(\bm{T}),z_{i+1}+c_{i+1}(\bm{T}))
\tph_{i+1}^{\bm{T}}\tph_i^{\bm{T}}
\nonumber\\
&\qquad=
\tph_{i+1}^{\bm{T}'}\tph_i^{\bm{T}'}
S_i(z_{i+1}+c_{i+1}(\bm{T}),z_i+c_i(\bm{T})).
\label{E2}
\end{align}
Here, the arguments of the $\tph$-factors on the right-hand side are obtained by making the above interchange of the $z$ and $w$ variables.

Using \eqref{E1} for the factors with indices greater than $i+1$, followed by \eqref{E2}, we obtain
\begin{align}
&S_i(z_i+c_i(\bm{T}),z_{i+1}+c_{i+1}(\bm{T}))
\Phi_{\bm{T}}(z_1,\ldots,z_n,w_1,\ldots,w_n)
\nonumber\\
&\qquad=
\left(
\prod_{k=1,\ldots,n}^{\leftarrow}
\tph_k^{\bm{T}'}(z'_1,\ldots,z'_k,w'_1,\ldots,w'_k)
\right)
\nonumber\\
&\qquad\qquad\times
S_i(z_{i+1}+c_{i+1}(\bm{T}),z_i+c_i(\bm{T}))
J_1(w_1+p_1(\bm{T}))\cdots J_n(w_n+p_n(\bm{T})).
\label{E3}
\end{align}

We now pass the Baxterized element on the right of \eqref{E3} through the $J$-factors. It commutes with all $J_k$ for $k\neq i,i+1$. For the remaining factors, we use
$$
s_iJ_i(a)=J_{i+1}(a)s_i,
\qquad
s_iJ_{i+1}(b)=J_i(b)s_i.
$$
Therefore,
\begin{align*}
&s_iJ_i(w_i+p_i(\bm{T}))
J_{i+1}(w_{i+1}+p_{i+1}(\bm{T}))\\
&\qquad=
J_i(w_{i+1}+p_{i+1}(\bm{T}))
J_{i+1}(w_i+p_i(\bm{T}))s_i.
\end{align*}

The scalar part of the Baxterized element is
$$
\frac{
\delta_{w_{i+1}+p_{i+1}(\bm{T}),
\,w_i+p_i(\bm{T})}
}{
z_{i+1}+c_{i+1}(\bm{T})-z_i-c_i(\bm{T})
}.
$$
When the Kronecker delta is zero, no scalar term occurs. When it is equal to one, we have
$$
w_{i+1}+p_{i+1}(\bm{T})
=
w_i+p_i(\bm{T}),
$$
and hence the product of the corresponding $J$-factors is unchanged by the interchange. Since $w'_i=w_{i+1}$ and $w'_{i+1}=w_i$, equation \eqref{E3} becomes
\begin{align}
&S_i(z_i+c_i(\bm{T}),z_{i+1}+c_{i+1}(\bm{T}))
\Phi_{\bm{T}}(z_1,\ldots,z_n,w_1,\ldots,w_n)
\nonumber\\
&\qquad=
\Phi_{\bm{T}'}(z'_1,\ldots,z'_n,w'_1,\ldots,w'_n)
S_i(z_{i+1}+c_{i+1}(\bm{T}),z_i+c_i(\bm{T})).
\label{E3-new}
\end{align}

The interchange of variables also gives
$$
(z'_1,\ldots,z'_n)\in H_{\bm{T}'}
\quad\Longleftrightarrow\quad
(z_1,\ldots,z_n)\in H_{\bm{T}}.
$$
If $i$ and $i+1$ occur in the same component, the standardness of both $\bm{T}$ and $\bm{T}'$ implies
$$
|c_i-c_{i+1}|\geq2.
$$
It follows from \eqref{Bx3} that the two Baxterized elements in \eqref{E3-new} are regular and invertible at $z_i=z_{i+1}$. If the entries occur in different components, their position parameters are different, and the Baxterized elements reduce to $s_i$. They are again invertible. Thus, $\Phi_{\bm{T}}$ is regular at $z_1=\cdots=z_n$ if and only if $\Phi_{\bm{T}'}$ is regular there.

The graph of standard bi-tableaux of a fixed shape, with edges given by standard adjacent transpositions, is connected. We can therefore pass from $\bm{T}$ to the hook bi-tableau $\bm{T}^{0}$ through a sequence of standard bi-tableaux. It remains to prove regularity for $\bm{T}^{0}$.

\medskip
\noindent
\textbf{Step 2: Removal of singularities for $\bm{T}^{0}$.}

Fix the hook bi-tableau $\bm{T}^{0}$ and let $c_j$ be the content of its entry $j$. Suppose that
$$
S_i(z_k+c_k,z_i+c_i)
$$
is a singularity term in $\tph_k$. This means that $i$ and $k$ lie on the same diagonal of one component of $\bm{T}^{0}$. In the ordered expression for $\tph_k$, this factor occurs immediately after
$$
S_{i+1}(z_k+c_k,z_{i+1}+c_{i+1}).
$$
The block $\tph_{i+1}$, which occurs later in the full product, contains
$$
S_i(z_{i+1}+c_{i+1},z_i+c_i).
$$

We bring this last factor next to the singularity term. For each
$\ell>i+1$, the braid and commuting relations allow us to pass it through $\tph_\ell$:
\begin{align*}
&\tph_\ell(z_1+c_1,\ldots,z_\ell+c_\ell,
w_1+p_1,\ldots,w_\ell+p_\ell)
S_i(z_{i+1}+c_{i+1},z_i+c_i)\\
&\qquad=
S_i(z_{i+1}+c_{i+1},z_i+c_i)
\tph'_\ell(z'_1+c'_1,\ldots,z'_\ell+c'_\ell,
w'_1+p'_1,\ldots,w'_\ell+p'_\ell),
\end{align*}
where $\tph'_\ell$ differs from $\tph_\ell$ only by the interchange of the data associated with $i$ and $i+1$.

Applying this relation successively for $
\ell=i+2,i+3,\ldots,k-1 $ places the factor immediately to the right of $\tph_k$. Passing it through the terminal word $s_1\cdots s_{k-1}$ changes its index and gives
$$
S_{i+1}(z_{i+1}+c_{i+1},z_i+c_i).
$$
The remaining commuting relations place it next to the singular factor. Hence, the part of the product containing the possible pole takes the form
\begin{align}
&S_{i+1}(z_k+c_k,z_{i+1}+c_{i+1})
S_i(z_k+c_k,z_i+c_i)
\nonumber\\
&\hspace{4cm}\times
S_{i+1}(z_{i+1}+c_{i+1},z_i+c_i).
\label{tri}
\end{align}

The special numbering of the hook bi-tableau now becomes important. The entries $i$ and $k$ are on the same diagonal, while $i+1$ is adjacent to $i$ in the row or column of the corresponding principal hook. Thus, $c_{i+1}=c_i-1$ in the column case, and $c_{i+1}=c_i+1$ in the row case.

On $H_{\bm{T}^{0}}$, in the column case,
$$
S_{i+1}(z_{i+1}+c_{i+1},z_i+c_i)^2 =
-2S_{i+1}(z_{i+1}+c_{i+1},z_i+c_i).
$$
Consequently, $-\frac12
S_{i+1}(z_{i+1}+c_{i+1},z_i+c_i)$ is an idempotent. In the row case,
$$
S_{i+1}(z_{i+1}+c_{i+1},z_i+c_i)^2 = 2S_{i+1}(z_{i+1}+c_{i+1},z_i+c_i),
$$
and hence $\frac12 S_{i+1}(z_{i+1}+c_{i+1},z_i+c_i)$ is an idempotent.

In the column case, the last factor of \eqref{tri} may therefore be written as
$$
-\frac12
S_{i+1}(z_{i+1}+c_{i+1},z_i+c_i)^2.
$$
Thus, one copy of this factor completes the triple \eqref{tri}, while the remaining idempotent factor can be moved back to its original position by reversing the preceding braid and commutation operations. The same reasoning applies in the row case, with $\frac12$ in place of $-\frac12$. Hence these rearrangements do not change the value of the full product.

The triple \eqref{tri} has the form considered in Lemma~\ref{Reg 1}. Indeed, on $H_{\bm{T}^{0}}$, the second and third spectral parameters differ by $\pm1$, while at $z_1=\cdots=z_n$ the equality $c_k=c_i$ makes the first and third parameters equal. Lemma~\ref{Reg 1} therefore shows that the possible pole in \eqref{tri} is removable.

Repeating this construction for every singular pair removes all the possible singularities in the product. Hence the restriction of $\Phi_{\bm{T}^{0}}$ to $H_{\bm{T}^{0}}$ is regular at $z_1=z_2=\cdots=z_n$. The reduction in Step 1 now proves the assertion for every standard bi-tableau $\bm{T}$.
\end{proof}

\begin{example}
Consider the hook bi-tableau $\bm{T}^0 = \left(\begin{ytableau}
       1& 3\\
       2& 4
\end{ytableau}, \quad \begin{ytableau}
       5& 6 
\end{ytableau} \right).$
The corresponding rational function is expressed as follows,
$$\Phi_{\bm{T}^0}(z_{1},\dots, z_{n}, w_{1}\dots, w_{n}) = \tph_6\tph_5\tph_4\tph_3\tph_2\tph_1 J_1 J_2 J_3 J_4 J_5 J_6.$$
For convenience, we write $\tph_k(z_{1}+c_{1},\dots, z_{k}+c_{k},w_{1}+p_{1}\dots, w_{k}+p_{k}) = \tph_k$ and $J_{k}(w_k+p_k) = J_k$. We notice that the singularity term only appears in the factor $\tph_4$, which is $S_1(z_4, z_1)$ in the product:
$$\Phi_{\bm{T}^0}(z_{1},\dots, z_{n}, w_{1}\dots, w_{n}) = \tph_6\tph_5S_3(z_4, z_3+1)S_2(z_4, z_2-1)S_1(z_4, z_1)s_1s_2s_3\tph_3\tph_2\tph_1 J_1 J_2 J_3 J_4 J_5 J_6.$$
Next pairing the singularity term with its corresponding triple term inside the bracket, we have, 
$$\Phi_{\bm{T}^0}(z_{1},\dots, z_{n}, w_{1}\dots, w_{n}) = \tph_6\tph_5.S_3(z_4, z_3+1)(S_2(z_4, z_2-1)S_1(z_4, z_1))s_1s_2s_3.\tph_3\tph_2\tph_1 J_1 J_2 J_3 J_4 J_5 J_6.$$
We then replace $S_2(z_4, z_2-1)S_1(z_4, z_1)$ in the product with $S_2(z_4, z_2-1)S_1(z_4, z_1)(-\frac{1}{2})S_2(z_2-1, z_1)$, where $(-\frac{1}{2})S_2(z_2-1, z_1)$ is an idempotent that does not alter the value of $\Phi_{\bm{T}^0}(z_{1},\dots, z_{n}, w_{1}\dots, w_{n})$. As the triple term is well-defined at $z_{1} = z_{2} = \dots = z_{n}$ and $w_i =0$ for $1\leq i\leq n$, the entire expression for $\Phi_{\bm{T}^0}(z_{1},\dots, z_{n}, w_{1}\dots, w_{n})$ remains well-defined.  
\end{example}

Therefore, based on the proposition above, we define an element $\Phi_{\bm{T}}\in \C[\Z_2\wr S_n]$ by evaluating \\$\Phi_{\bm{T}}(z_{1},\dots, z_{n}, w_{1}\dots, w_{n})$ at $z_1 = z_2 =\dots = z_n$ and $w_1 = w_2 =\dots = w_n = 0$. To prove the next part of Theorem 2.2, we will need the following propositions.

\begin{proposition}\label{prop:unit-coefficient}
The coefficient of $\Phi_{\bm{T}}\in \C[\Z_2\wr S_n]$ at the unit element of $\C[\Z_2\wr S_n]$ is $\frac{1}{2^{n}}$.
\end{proposition}

\begin{proof}
We have the expression: $$\Phi_{\bm{T}} = \tph_{n}\tph_{n-1}\dots\tph_{1}J_{1}(p_{1})J_{2}(p_2)\dots J_{n}(p_{n}),$$
where \begin{align*}\tph_{k} = &S_{k-1}(c_k, c_{k-1})\dots S_{1}(c_k, c_1)s_{1}\dots s_{k-1}\\
 = &(s_{k-1} + \frac{\delta_{p_{k},p_{k-1}}}{c_{k}-c_{k-1}})\dots (s_{1} + \frac{\delta_{p_{k},p_{1}}}{c_{k}-c_{1}})s_{1}\dots s_{k-1}.
\end{align*}
The coefficient of the product $(s_{k-1} + \frac{\delta_{p_{k},p_{k-1}}}{c_{k}-c_{k-1}})\dots (s_{1} + \frac{\delta_{p_{k},p_{1}}}{c_{k}-c_{1}})$ at $s_{k-1}\dots s_{1}$ is 1. This implies that the coefficient of $\tph_k$ at the unit element is 1. Therefore, the coefficient of $\tph_{n}\tph_{n-1}\dots\tph_{1}$ at unit element is 1. However, the coefficient of the product $J_{1}(p_{1})J_{2}(p_2)\dots J_{n}(p_{n})$ at unit element is $\frac{1}{2^{n}}.$ Hence, the coefficient of $\Phi_{\bm{T}}$ at unit element of $\Z_{2}\wr S_{n}$ is $\frac{1}{2^{n}}$. 
\end{proof}

Note that $\frac{1}{2^{n}} = \frac{1}{f_{\bm{\lambda}}}\frac{d_{\bm{\lambda}}}{|\Z_2\wr S_n)|},$ where $d_{\bm{\lambda}}$ is the dimension of the irreducible representation corresponding to the bi-partition $\bm{\lambda}$.

\begin{proposition}\label{prop:divisibility}
Suppose that $i$ and $i+1$ are adjacent in the standard bi-tableau $\bm T$. If they lie in the same row, then $\Phi_{\bm T}$ is divisible on the left by $s_i+1$. If they lie in the same column, then $\Phi_{\bm T}$ is divisible on the left by $s_i-1$.
\end{proposition}

\begin{proof}
In this proof, all rational functions are considered after their restriction to the corresponding hook subspace $H_{\bm T}$.

Using the factorization in (\ref{RF3}), we have
$$
\Phi_{\bm T}(\bm z,\bm w)
=
\left(
\prod_{k=1,\ldots,n}^{\leftarrow}
\tph_k\bigl(z_1+c_1,\ldots,z_k+c_k,
w_1+p_1,\ldots,w_k+p_k\bigr)
\right)
J_1(w_1+p_1)\cdots J_n(w_n+p_n).
$$
Since the proposition concerns the final value obtained after setting
$w_1=\cdots=w_n=0$, define
$$
\Psi_{\bm T}(\bm z)
:=
\prod_{k=1,\ldots,n}^{\leftarrow}
\tph_k\bigl(z_1+c_1,\ldots,z_k+c_k,
p_1,\ldots,p_k\bigr).
$$
It follows that
$$
\left.
\Phi_{\bm T}(\bm z,\bm w)
\right|_{w_1=\cdots=w_n=0}
=
\Psi_{\bm T}(\bm z)
J_1(p_1)\cdots J_n(p_n).
$$
The $J$-factors appear on the right. Therefore, it is sufficient to establish the required left divisibility for the evaluated value of $\Psi_{\bm T}(\bm z)$. Set
$$
x_j=z_j+c_j,\qquad 1\leq j\leq n.
$$
By the definition of $\tph_k$,
$$
\tph_k
=
S_{k-1}(x_k,x_{k-1};p_k,p_{k-1})
\cdots
S_1(x_k,x_1;p_k,p_1)
s_1\cdots s_{k-1}.
$$
Moving the generators successively through the terminal word
$s_1\cdots s_{k-1}$ gives
$$
\begin{aligned}
\tph_k
&=
\left(
1+\frac{\delta_{p_k,p_{k-1}}(k-1\,\,k)}
{x_k-x_{k-1}}
\right)
\cdots
\left(
1+\frac{\delta_{p_k,p_1}(1\,\,k)}
{x_k-x_1}
\right)\\
&=
\left(
1-\frac{\delta_{p_k,p_{k-1}}(k-1\,\,k)}
{x_{k-1}-x_k}
\right)
\cdots
\left(
1-\frac{\delta_{p_k,p_1}(1\,\,k)}
{x_1-x_k}
\right).
\end{aligned}
$$
Consequently, every nonidentity factor in $\Psi_{\bm T}(\bm z)$ is of the form
$$
1-\frac{(l\,\,k)}{x_l-x_k},
\qquad 1\leq l<k\leq n,
$$
and it occurs only when $p_l=p_k$. For $a\in\{1,2\}$, let
$$
I_a=\{m_1^{(a)}<\cdots<m_{n_a}^{(a)}\}
$$
be the set of entries contained in the $a$-th component of $\bm T$. Since the two components have distinct position parameters,
$$
p_l=p_k
\quad\Longleftrightarrow\quad
l\text{ and }k\text{ occur in the same component of }\bm T.
$$
Thus, a factor involving entries from different components is equal to $1$. For each $a\in\{1,2\}$, denote by
$\Psi_{\bm T}^{(a)}(\bm z)$ the ordered product of all factors
$$
1-\frac{(l\,\,k)}{x_l-x_k}
$$
for which $l,k\in I_a$. The two products have disjoint supports and hence commute. We may therefore write
$$
\Psi_{\bm T}(\bm z)
=
\Psi_{\bm T}^{(1)}(\bm z)
\Psi_{\bm T}^{(2)}(\bm z).
$$

We next interpret each of these two products as a symmetric-group fusion function. Fix $a\in\{1,2\}$, and let $U^{(a)}$ be the ordinary standard tableau obtained from the $a$-th component of $\bm T$ by replacing
$$
m_1^{(a)}<\cdots<m_{n_a}^{(a)}
$$
with $1<\cdots<n_a$, while keeping their boxes unchanged. Then
$$
c_j\bigl(U^{(a)}\bigr)
=
c_{m_j^{(a)}}(\bm T),
\qquad 1\leq j\leq n_a.
$$
Set $u_j=z_{m_j^{(a)}}$ for $1\leq j\leq n_a$. It follows that $x_{m_j^{(a)}}=
u_j+c_j\bigl(U^{(a)}\bigr)$. Consider the algebra embedding
$$
\iota_a:\C[S_{n_a}]\longrightarrow\C[S_n]
$$
defined by
$$
\iota_a((j\,\,l))
=
(m_j^{(a)}\,\,m_l^{(a)}).
$$
We extend $\iota_a$ coefficientwise to rational functions. The factor in
$\Psi_{\bm T}^{(a)}(\bm z)$ corresponding to
$m_j^{(a)}<m_l^{(a)}$ can then be expressed as
\begin{align*}
1-\frac{(m_j^{(a)}\,\,m_l^{(a)})}
{x_{m_j^{(a)}}-x_{m_l^{(a)}}}
=
\iota_a\left(
1-\frac{(j\,\,l)}
{u_j+c_j(U^{(a)})-u_l-c_l(U^{(a)})}
\right).
\end{align*}
The Yang--Baxter and commuting relations allow us to arrange these factors in the usual lexicographic order. Hence
$$
\Psi_{\bm T}^{(a)}(\bm z)
=
\iota_a\left(
F_{U^{(a)}}(u_1,\ldots,u_{n_a})
\right),
$$
where
$$
F_{U^{(a)}}(u_1,\ldots,u_{n_a})
=
\prod_{1\leq j<l\leq n_a}^{\rightarrow}
\left(
1-\frac{(j\,\,l)}
{u_j+c_j(U^{(a)})-u_l-c_l(U^{(a)})}
\right)
$$
is the symmetric-group fusion function associated with $U^{(a)}$.

By the hook fusion procedure for the symmetric group, the restriction of this function to the hook subspace associated with $U^{(a)}$ is regular when all the variables are equal; see
\cite[Theorem~1.2]{Gri05}. Its evaluated value is divisible on the left by $(1+(q\,\,q+1))$ when $q$ and $q+1$ are adjacent in the same row, and by $(1-(q\,\,q+1))$ when they are adjacent in the same column; see
\cite[Corollary~3.13]{Gri07}.

Let the entries $i$ and $i+1$ lie in the $a$-th component of $\bm T$. Since they are consecutive integers, they occur consecutively in the increasing sequence
$$
m_1^{(a)}<\cdots<m_{n_a}^{(a)}.
$$
Thus, for some $q$, $m_q^{(a)}=i$ and $m_{q+1}^{(a)}=i+1$. In particular,$\iota_a((q\,\,q+1))=(i\,\,i+1)=s_i$. Suppose first that $i$ and $i+1$ are adjacent in the same row of $\bm T$. Then $q$ and $q+1$ are adjacent in the same row of $U^{(a)}$. Therefore, after evaluation,
$$
\left.
F_{U^{(a)}}\right|_{u_1=\cdots=u_{n_a}}
=
\bigl(1+(q\,\,q+1)\bigr)B_a
$$
for some $B_a\in\C[S_{n_a}]$. Applying $\iota_a$, we obtain
$$
\left.
\Psi_{\bm T}^{(a)}(\bm z)
\right|_{z_1=\cdots=z_n}
=
(1+s_i)\iota_a(B_a).
$$

Let $b$ be the other component. Since
$\Psi_{\bm T}^{(b)}(\bm z)$ is supported only on the entries outside the $a$-th component, its evaluated value commutes with $s_i$. Consequently, for some $A\in\C[S_n]$,
\begin{align*}
\left.
\Psi_{\bm T}(\bm z)
\right|_{z_1=\cdots=z_n}
&=
\left.
\Psi_{\bm T}^{(a)}(\bm z)
\right|_{z_1=\cdots=z_n}
\left.
\Psi_{\bm T}^{(b)}(\bm z)
\right|_{z_1=\cdots=z_n}\\
&=
(1+s_i)A.
\end{align*}
It now follows that
\begin{align*}
\Phi_{\bm T}
&=
\left.
\Phi_{\bm T}(\bm z,\bm w)
\right|_{\substack{z_1=\cdots=z_n\\
w_1=\cdots=w_n=0}}\\
&=
\left.
\Psi_{\bm T}(\bm z)
\right|_{z_1=\cdots=z_n}
J_1(p_1)\cdots J_n(p_n)\\
&=
(1+s_i)
\bigl(AJ_1(p_1)\cdots J_n(p_n)\bigr).
\end{align*}
Therefore, $\Phi_{\bm T}$ is divisible on the left by $s_i+1$.

Now suppose that $i$ and $i+1$ are adjacent in the same column of $\bm T$. In this case, $q$ and $q+1$ are adjacent in the same column of $U^{(a)}$. By
\cite[Corollary~3.13]{Gri07}, the evaluated symmetric-group fusion function associated with $U^{(a)}$ is divisible on the left by $
(1-(q\,\,q+1))$. After applying $\iota_a$ and including the factor belonging to the other component, we find that
$$
\left.
\Psi_{\bm T}(\bm z)
\right|_{z_1=\cdots=z_n}
$$
is divisible on the left by $1-s_i$. Since $
1-s_i=-(s_i-1)$, it is also divisible on the left by $s_i-1$. Multiplying on the right by $J_1(p_1)\cdots J_n(p_n)$ does not affect the left divisibility. Hence, $\Phi_{\bm T}$ is divisible on the left by $s_i-1$.
\end{proof}

We now prove the second statement of the Theorem 2.2.

\begin{proposition}\label{prop:final-evaluation}
The value of the restriction to $H_{\bm{T}}$ of the rational function $\Phi_{\bm{T}}(z_{1},\dots, z_{n}, w_{1}\dots, w_{n})$  at $z_{1} = z_{2} = \dots = z_{n}$ and $w_{1} = w_{2} = \dots = w_{n} = 0$ coincides with $\frac{1}{2^n}F_{\bm{T}}$.
\end{proposition}

\begin{proof}
We compare the coefficient of every element of $\Z_2\wr S_n$ in the evaluated function $\Phi_{\bm T}$ with the corresponding coefficient in the diagonal matrix element $F_{\bm T}$. Write
$$
\Phi_{\bm T}
=
\sum_{g\in\Z_2\wr S_n}\alpha_g(\bm T)\,g
$$
and
$$
F_{\bm T}
=
\sum_{g\in\Z_2\wr S_n}
\langle v_{\bm T},gv_{\bm T}\rangle\,g.
$$

\medskip
\noindent
\textbf{Step 1: The normalization.}

By Proposition~\ref{prop:unit-coefficient}, the coefficient of the identity element in $\Phi_{\bm T}$ is
$$
\alpha_1(\bm T)=\frac{1}{2^n}.
$$
On the other hand,
$$
\langle v_{\bm T},v_{\bm T}\rangle=1,
$$
so the identity coefficient in $F_{\bm T}$ is $1$. Thus, the required normalization factor is necessarily $\frac{1}{2^n}$.

\medskip
\noindent
\textbf{Step 2: The generators of the normal abelian subgroup.}

The evaluated function has the form
$$
\Phi_{\bm T}
=
\tph_n\cdots\tph_1
J_1(p_1)\cdots J_n(p_n),
$$
where $\tph_k\in\C[S_n]$ and
$$
J_i(p_i)=\frac12(1+p_ij_i).
$$
The commuting elements $j_1,\ldots,j_n$ generate a subgroup isomorphic to $(\Z_2)^n$. Every element of this subgroup has a unique expression
$$
j_1^{\varepsilon_1}\cdots j_n^{\varepsilon_n},
\qquad
\varepsilon_i\in\{0,1\}.
$$
Moreover, every element of $\Z_2\wr S_n$ can be written uniquely as
$$
\sigma j_1^{\varepsilon_1}\cdots j_n^{\varepsilon_n},
\qquad
\sigma\in S_n.
$$

Since $\tph_n\cdots\tph_1$ belongs to $\C[S_n]$, the terms involving the elements $j_i$ come from the product $J_1(p_1)\cdots J_n(p_n)$. To obtain $j_i$, we must choose the term $p_ij_i$ from $J_i(p_i)$ and the identity term from every other $J_k(p_k)$. Therefore,
$$
\alpha_{j_i}(\bm T)=\frac{1}{2^n}p_i.
$$

Since $p_i\in\{1,-1\}$ and $j_i^2=1$, we have $J_i(p_i)j_i=p_iJ_i(p_i)$. It follows that
\begin{equation}
\Phi_{\bm T}j_i=p_i\Phi_{\bm T}.
\label{CR1}
\end{equation}
The seminormal action also gives $j_iv_{\bm T}=p_iv_{\bm T}$. Consequently, the diagonal matrix element satisfies
$$
F_{\bm T}j_i=p_iF_{\bm T}.
$$
In particular, the coefficient of $t=j_1$ in $F_{\bm T}$ is $p_1$, whereas its coefficient in $\Phi_{\bm T}$ is $\frac{1}{2^n}p_1$.

More generally, relation \eqref{CR1} implies
$$
\alpha_{\sigma j_1^{\varepsilon_1}\cdots
j_n^{\varepsilon_n}}(\bm T) = p_1^{\varepsilon_1}\cdots p_n^{\varepsilon_n}
\alpha_\sigma(\bm T).
$$
The coefficients of $F_{\bm T}$ satisfy the same identity. Therefore, it remains to compare the coefficients corresponding to permutations $\sigma\in S_n$.

\medskip
\noindent
\textbf{Step 3: The symmetric-group generators.}

We consider separately the two possible positions of $i$ and $i+1$ in $\bm T$.

\smallskip
\noindent
\emph{Case (a): The entries $i$ and $i+1$ lie in the same component.}

Suppose first that $s_i\bm T$ is non-standard. Then $i$ and $i+1$ are adjacent in the same row or in the same column. By Proposition~\ref{prop:divisibility},
\begin{equation}
\Phi_{\bm T}=\pm s_i\Phi_{\bm T},
\label{CR2}
\end{equation}
where the positive sign occurs in the row case and the negative sign occurs in the column case. Comparing coefficients in \eqref{CR2}, we obtain $$
\alpha_{s_ig}(\bm T) =\pm\alpha_g(\bm T)
$$
for every $g\in\Z_2\wr S_n$. These are precisely the relations satisfied by the coefficients of $F_{\bm T}$ in the corresponding row and column cases.

Now suppose that $\bm T'=s_i\bm T$ is standard. Evaluating relation \eqref{E3-new} at $z_1=\cdots=z_n$ and $w_1=\cdots=w_n=0$, gives
$$
S_i(c_i,c_{i+1})\Phi_{\bm T}=\Phi_{\bm T'}S_i(c_{i+1},c_i).
$$
Equivalently,
\begin{align*}
(1-d_i(\bm T')^2)\Phi_{\bm T}=(s_i-d_i(\bm T'))\Phi_{\bm T'}(s_i-d_i(\bm T')),
\end{align*}
where $d_i(\bm T)=(c_{i+1}-c_i)^{-1}$. The coefficient of $g$ on both sides gives
\begin{align}
(1-d_i(\bm T')^2)\alpha_g(\bm T)=\alpha_{s_igs_i}(\bm T')-d_i(\bm T')\alpha_{s_ig}(\bm T')-d_i(\bm T')\alpha_{gs_i}(\bm T')+d_i(\bm T')^2\alpha_g(\bm T').
\label{CR3}
\end{align}
Taking the identity coefficient into account, we obtain $\alpha_{s_i}(\bm T)=
\frac{1}{2^n}d_i(\bm T)$. The seminormal action gives $\langle v_{\bm T},s_iv_{\bm T}\rangle=d_i(\bm T)$. Thus, the coefficients of $s_i$ in $\Phi_{\bm T}$ and $F_{\bm T}$ differ by the factor $\frac{1}{2^n}$.

\smallskip
\noindent
\emph{Case (b): The entries $i$ and $i+1$ lie in different components.}

In this case, $\bm T'=s_i\bm T$ is standard and the Baxterized element reduces to $s_i$. Therefore, relation \eqref{E3-new} becomes
\begin{equation}
\Phi_{\bm T}
=
s_i\Phi_{\bm T'}s_i.
\label{CR4}
\end{equation}

The componentwise decomposition obtained in the proof of Proposition~\ref{prop:divisibility} shows that the symmetric part of $\Phi_{\bm T}$ contains only those permutations that preserve the sets of entries belonging to the two components. However, the transposition $s_i=(i\,\,i+1)$ interchanges entries from different components. Hence it cannot occur in the symmetric part of $\Phi_{\bm T}$. By the uniqueness of the decomposition $\Z_2\wr S_n=(\Z_2)^n\rtimes S_n$, no term containing a nonidentity element of $(\Z_2)^n$ can contribute to the coefficient of $s_i$. Therefore, $\alpha_{s_i}(\bm T)=0$. 

The seminormal action in this case is
$s_iv_{\bm T}=v_{\bm T'}$, and the orthogonality of the seminormal basis gives
$\langle v_{\bm T},s_iv_{\bm T}\rangle =
\langle v_{\bm T},v_{\bm T'}\rangle =0$.
Thus, the coefficients again agree up to the factor $\frac{1}{2^n}$.

Relations \eqref{CR1}--\eqref{CR4} give the same recursive relations as those satisfied by the coefficients of the diagonal matrix elements. The graph of standard bi-tableaux of a fixed shape, whose edges correspond to standard adjacent transpositions, is connected. Starting from any fixed standard bi-tableau and applying these relations along such a sequence determines the coefficients for every standard bi-tableau of the same shape. Their compatibility with the Coxeter relations ensures that the result is independent of the chosen sequence.

The coefficients of $\Phi_{\bm T}$ and $\frac{1}{2^n}F_{\bm T}$ satisfy the same relations and have the same identity coefficient. Hence they coincide for every element of $\Z_2\wr S_n$. Therefore,
$$
\Phi_{\bm T}
=
\frac{1}{2^n}F_{\bm T}.
$$
\end{proof}

\noindent \textbf{Acknowledgements}

The first author's research is supported by Indian Institute of Science Education and Research Thiruvananthapuram PhD fellowship. The second author's research was partially supported by IISER-Thiruvananthapuram, SERB-Power Grant SPG/2021/004200.


\begin{thebibliography}{nyt}
\bibliographystyle{abbrv}

\bibitem{Ar97}
Ram, Arun, Seminormal representations of Weyl groups and Iwahori-Hecke algebras, {\it Proc. London Math. Soc. (3)} {\bf 75} (1997) 99--133.

\bibitem{Ca98}
Can, H., On the inequivalence and standard basis of the Specht modules of the hyperoctahedral groups, {\it Commun. Fac. Sci. Univ. Ank. Ser. A1 Math. Stat.} {\bf 47} (1998) 125--138.

\bibitem{Che86}
Cherdnik, I.V., Special bases of irreducible representations of a degenerate affine Hecke algebra, {\it Funktsional. Anal. i Prilozhen.} {\bf 20} (1986) 87--88.

\bibitem{Gri05}
Grime, James, The hook fusion procedure, {\it Electron. J. Combin.} {\bf 12} (2005) 26, 14.

\bibitem{Gri07}
Grime, James, The hook fusion procedure for Hecke algebras, {\it Journal of algebra} {\bf 309} (2005) 744--759.

\bibitem{Ho74}
Hoefsmit, Peter Norbert, Representations of Hecke algebras of finite groups with BN-pairs of classical type, {\it Thesis (Ph.D.)--The University of British Columbia (Canada)} (1974).

\bibitem{IM08}
Isaev, A.P. and Molev, A.I. and Os'Kin, A.F., On the idempotents of Hecke algebras, {Lett. Math. Phys.} {\bf 85} (2008) 79--90.

\bibitem{IM10}
Isaev, A.P. and Molev, A.I., Fusion procedure for the Brauer algebra, {Algebra i Analiz} {\bf 22} (2010) 142--154.

\bibitem{IM14}
Isaev, A.P. and Molev, and Ogievetsky, O. V., Idempotents for Birman-Murakami-Wenzl algebras and reflection equation, {Adv. Theor. Math. Phys.} {\bf 18 } (2014) 1--25.

\bibitem{Ju66}
Jucys, A.A., On the Young operators of symmetric groups, {Litovsk. Fiz. Sb.} {\bf 6 } (1966) 163--180.

\bibitem{Mo08}
Molev, A.I., On the fusion procedure for the symmetric group, {Rep. Math. Phys.} {\bf 61 } (2008) 181--188.

\bibitem{Naz97}
Nazarov, Maxim, Young's symmetrizers for projective representations of the symmetric group, {Adv. Math.} {\bf 127 } (1997) 190--257.


\bibitem{Naz98}
Nazarov, Maxim, Yangians and Capelli identities, {Kirillov's seminar on representation theory, Amer. Math. Soc. Transl. Ser. 2} {\bf 181 } (1998) 139--163.

\bibitem{Naz04}
Nazarov, Maxim, Representations of twisted Yangians associated with skew Young diagrams, {Selecta Math. (N.S.)} {\bf 10 } (2004) 71--129.

\bibitem{NT02}
Nazarov, Maxim and Tarasov, Vitaly, On irreducibility of tensor products of Yangian modules associated with skew Young diagrams, {Duke Math. J.} {\bf 112 } (2002) 343--378.

\bibitem{OP14}
Ogievetsky, O. V. and Poulain d'Andecy, L., Fusion procedure for Coxeter groups of type B and complex reflection groups {$G(m,1,n)$}, {Proc. Amer. Math. Soc.} {\bf 142 } (2014) 2929--2941.

\bibitem{PL17}
Poulain d'Andecy, L., Fusion formulas and fusion procedure for the Yang-Baxter equation, {Algebr. Represent. Theory} {\bf 20 } (2017) 1379--1414.

\bibitem{Pu99}
Pushkarev, Igor A., On the representation theory of wreath products of finite groups and symmetric groups, {Journal of Mathematical Sciences} {\bf 96 } (2017) 3590--3599.
\end{thebibliography}
\end{document}